\documentclass[twoside,12pt]{article}
\usepackage{amsmath,amssymb,epsfig}
\date{}
\newtheorem{Lemma}{LEMMA}[section]
\newtheorem{Corollary}[Lemma]{COROLLARY}
\newtheorem{Theorem}[Lemma]{THEOREM}
\newtheorem{Proposition}[Lemma]{PROPOSITION}
\newtheorem{Definition}[Lemma]{DEFINITION}
\newtheorem{Example}[Lemma]{Example}

\newcommand{\bnum}{\begin{enumerate}}
\newcommand{\enum}{\end{enumerate}}
\newcommand{\bi}{\begin{itemize}}
\newcommand{\ei}{\end{itemize}}
\newcommand{\btab}{\begin{tabular}}
\newcommand{\etab}{\end{tabular}}
\newcommand{\beq}{\begin{eqnarray*}}
\newcommand{\eeq}{\end{eqnarray*}}
\newcommand{\beqn}{\begin{eqnarray}}
\newcommand{\eeqn}{\end{eqnarray}}

\newcommand{\bq}{\begin{equation}}
\newcommand{\eq}{\end{equation}}

\newcommand{\CC}{{\cal C}}
\newcommand{\CU}{{\cal U}}

\newcommand{\CL}{{\cal L}}
\newcommand{\CM}{{\cal M}}

\newcommand{\CP}{{\cal P}}

\newcommand{\CS}{{\cal S}}

\def\phi{\varphi}
\def\epsilon{\varepsilon}
\newcommand{\BS}{\mathbb S}
\newcommand{\BR}{\mathbb R}
\newcommand{\BC}{\mathbb C}
\newcommand{\BF}{\mathbb F}

\newcommand{\BH}{\mathbb H}
\newcommand{\BO}{\mathbb O}
\newcommand{\kasten}{\vbox{\hrule height 8pt width 8.6pt depth -7.4pt
    \hbox{\vrule width 0.6pt height 7.4pt
    \kern 7.4pt \vrule width 0.6pt height 7.4pt}
    \hrule height 0.6pt width 8.6pt}}
\newcommand{\ok}{\hfill\kasten}
\newcommand{\bpf}{\begin{Proof}}
\newcommand{\epf}{\ok\end{Proof}\bigskip\noindent}
\newcommand{\bthm}{\begin{Theorem}}
\newcommand{\ethm}{\end{Theorem}}
\newcommand{\ble}{\begin{Lemma}}
\newcommand{\ele}{\end{Lemma}}
\newcommand{\bprop}{\begin{Proposition}}
\newcommand{\eprop}{\end{Proposition}}
\newcommand{\bcor}{\begin{Corollary}}
\newcommand{\ecor}{\end{Corollary}}

\begin{document}
\title{Locally classical stable planes}
\author{Rainer L\"owen }

\maketitle
\thispagestyle{empty}

\begin{abstract}

We show that a simply connected stable plane with connected lines 
is isomorphic to an open subplane of a 
classical projective plane (i.e., a plane over the real or complex numbers, the 
quaternions or the octonions) if it has that property locally. We also give several examples
indicating that our hypotheses cannot be relaxed much further.

\end{abstract}

\section{Introduction}\label{intro}

The classical compact projective planes are the planes $\CP_2\BF$ over the algebras $\BF \in
\{\BR, \BC, \BH, \BO\}$ (real and complex numbers, quaternions and octonions).
Our aim is to characterize open subplanes of these projective planes among the 
general class of stable planes  by a local condition.
See Section \ref{prelim} for  definitions.\

Writing this paper was motivated by several observations. Highly sophisticated local 
conditions characterizing
the classical planes among differentiable planes were studied by Patrick Leymann 
in his recent doctoral 
dissertation \cite{leymann}. This reminded me of the much simpler, older results by C. Polley 
\cite{polley-szeb},  \cite{polley-lok}, \cite{polley-moeb} characterizing 
open subplanes of $\CP_2\BR$ 
among 2-dimensional stable planes with connected lines
using either a local Desargues condition or 
local isomorphism to a 
classical plane. Originally, my aim was merely to extend Polley's results 
to higher dimensions. Regarding the local isomorphism condition, this
is made possible by two new ingredients. The first one is
the `local fundamental theorem' \cite{locfund}, which replaces H. Salzmann's result 
\cite {sz67}, 1.1 used by Polley in the 2-dimensional case. The second one is
the isotopy theorem for arcs 
by J. Martin and D. Rolfsen \cite{martin}, which is needed here 
for a monodromy type argument. 
It turns out that the connectedness condition for 
lines can be relaxed; it suffices to require that the connected lines form a dense subset of 
the line space.
The Desargues condition entered into Polley's 
setup via a result by Busemann \cite{busemann},
which works only when lines can be considered as geodesics, i.e., 
in planes of dimension 2. 
In higher dimensions, this would require different methods, 
and will not be considered here.  \

While checking for possible other obstacles, I noticed that Polley's proof is 
written rather sketchily in places, and that there is a more serious issue concerning 
applications of his Lemma 1 together with Folgerung 1, which correspond to our 
Lemma \ref{two}, assertions (a) and (b), respectively. 
The latter assertion applies to lines meeting the 
intersection of two given 
open sets, but in \cite{polley-szeb} it is used several times 
in situations where this is not known. 
We circumvent this problem by introducing \it local 
homomorphisms \rm of planes and studying 
their properties, see Section \ref{hom}. We never use assertion (b), substituting it 
with Lemma \ref{zush} on properties of local homomorphisms.
The crucial assertion (a) allows us to adapt two homomorphisms with overlapping open 
domains so as to create a homomorphism on the union of the domains. 
It is obtained in higher dimensions using the 
local fundamental theorem, which generalizes the result by H. Salzmann 
used in \cite{polley-szeb}.
\

Finally, the 
existence of a non-classical but locally classical plane on a cylinder surface convinced me 
that a thorough revision of the original proof was indicated 
even in dimension 2. 
That plane is an open subplane of a 
Moulton projective plane, see Example \ref{excylinder}. It has some disconnected lines, but the 
connected ones are dense in the line space. The impression gained here
is that the 2-dimensional case is the most complicated one. A specific reason for this is the 
existence of the Moulton planes \cite{moulton}, see also \cite{CPP}, Section 34,  
which are responsible for almost all the examples 
given in Section \ref{ex}. There is no analogue of the Moulton planes in higher dimensions, 
see \cite{CPP}, 64.18. The reason why the proofs fail if the point set is a cylinder is the 
existence of pairs of arcs in the cylinder that are homotopic, but not isotopic, relative end 
points \cite{feustel}. This means in particular that special care is needed with all 
arguments that rely on isotopies of arcs. This problem is non-existent in higher 
dimensions by the result of Martin and Rolfsen \cite{martin} mentioned above.\\

\section{Preliminaries}\label{prelim}

A \it stable plane \rm is a pair $\CS = (S,\CL)$ consisting of a locally compact 
topological space $S$ of positive topological dimension (the \it point set\rm ) 
and a set $\CL$ of closed
subsets of $S$, called \it lines\rm, such that every pair  of distinct points $x,y\in S$
is contained in a unique line $x \vee y \in \CL$. 
The line set $\CL$ is supposed to carry a topology such that the join operation $\vee$ is 
continuous, and such that the set of intersecting line pairs is open in $\CL \times \CL$ 
(this is referred to as the \it stability condition\rm ) and moreover,
the intersection operation $\wedge$ defined on this set of pairs is continuous.
In a stable plane, a compact line intersects every other line \cite{diss}, 1.15; 
hence $\CS$ is a compact projective plane if and only if all its lines are compact. \

With the exception of the extension theorem \ref{ext}, we shall exclusively deal with planes 
that are locally isomorphic to a classical projective plane $\CP_2\BF$ 
as defined in Section \ref{hom}.
This implies that lines are manifolds of dimension $l \in \{1,2,4,8\}$ and the point set 
is a manifold of dimension $2l$. 
The set $\CL_x$ of all lines passing through a point $x$
is then homeomorphic to the sphere $\BF \cup {\infty} \approx \BS_l$. \\

\section{Homomorphisms and local homomorphisms}\label{hom}

\bf DEFINITION \rm  Let $\CS=(S,\CL)$ and $\CS'=(S',\CL')$ be stable planes 
of the same dimension. \

(a) A \it homomorphism  \rm $\CS \to \CS'$ is a pair of continuous maps  
$(\phi,\Phi): (S,\CL)\to (S',\CL')$ such that for each line $L \in \CL$ the 
image $L^\phi = \{x^\phi \ \vert \ x \in L\}$ is contained in the line $L^\Phi$; 
moreover, we require that $\phi$ is locally injective.
The map $\Phi$ is of course determined by $\phi$ and will often be suppressed in the notation.
Also continuity of $\Phi$ follows easily from that of $\phi$. \

(b) A map $\phi: S \to S'$ is called a \it local homomorphism \rm if every point $p\in S$ 
has an open neighborhood $U_p$ in $S$ such that $\phi$ defines a homomorphism $\CU_p \to \CS'$,
where $\CU_p$ is the plane induced on $U_p$ by $\CS$. 
\\

It is easy to produce examples of disconnected planes and local homomorphisms 
defined on them that are not homomorphisms.
See Example \ref{exdisconn} for such an example defined on a connected plane. 
We stress the point that in this plane the connected lines do not form a dense set.
We show next that a homomorphism is in fact globally injective. It seems 
unlikely that this holds for local homomorphisms, but I do not know a counterexample.\\

\ble\label{inj}
Every  homomorphism of stable planes is injective.
\ele

In other words, a homomorphism $\CS \to \CS'$ can be seen as an isomorphism of $\CS$
onto the plane induced by $\CS'$ on the set $S^\phi \subseteq S'$, which is 
open by domain invariance. 

\bpf If not, choose distinct points $x,y,z \in S$ not on one line, 
such that $x^\phi = y^\phi \ne z^\phi$. 
Let $L = z \vee x$ and $M = z \vee y$; then $L \ne M$ but $L^\Phi = M^\Phi$. A point $p \notin L \cup M$ 
is contained in a line meeting both $L$ and $M$ in points $u,v$, respectively, 
which are close to $z$
and have distinct images by local injectivity. Then $p^\phi \in u^\phi \vee v^\phi = M^\Phi$,
so $S^\phi \subseteq M^\Phi$ contrary to the definition of a homomorphism.
\epf

\begin{Definition} \rm
Finite-dimensional stable planes exist in dimensions 2,4,8, and 16 only \cite{crelle}. 
We call a stable plane \it classical \rm if it is isomorphic to an open subplane of 
the classical projective plane $\CP$ of the same dimension, or in other words, if it 
admits a homomorphism into $\CP$. \

A stable plane $\CS$ is said to be 
\it locally classical \rm if every point admits a neighborhood which defines a 
classical subplane $\CU$ of $\CS$. Again this is equivalent to the existence of a 
homomorphism $\CU \to \CP$.
\end{Definition}

It is easy to see that a homomorphism is determined by its restriction to any open subplane.
This is obviously not true for local homomorphisms: Let $\CS$ be a subplane of the real 
affine plane whose point set consists of two disjoint disks; map one of them into the 
affine plane by the identity map and rotate the other. Yet we can prove the following.

\ble\label{extuniq}
Let $\phi, \psi : \CS \to \CS'$ be two local homomorphisms that agree on some open 
set $U \subseteq S$. If $S$ is connected, then it follows that $\phi = \psi$. 
\ele

\bpf
Let $X =\{x \in S \ \vert \ x^\phi = x^\psi\}$. The interior $Y = \mathop{\mathrm{int}} Y$ is 
nonempty by assumption. It suffices to show that the closure of $Y$ is contained in $Y$.  
So let $p$ be a boundary point of $Y$. There is an open neighborhood $W$ of $p$ such 
that both $\phi$ and $\psi$ induce homomorphisms on $W$. Now $W$ contains an open subset 
$Z$ of $X$, and $p$ can be obtained as the intersection of two lines meeting $Z$. 
Both of these lines have the same image under $\phi$ and $\psi$, 
hence $p^\phi = p^\psi$.
\epf

The following  property of local homomorphisms is the key to our proof of Theorem \ref{main}.

\ble\label{zush}
The image of a connected line $L$ under a local homomorphism is contained in a single line.
\ele

\bpf
Every point of $L$ has a neighborhood whose image is contained in a line. 
Let $M$ be one of those lines. Then the set $X$ of points of $L$ that are
mapped into $M$ has a nonempty interior $Y$ in $L$.
We claim that $Y$ has empty boundary in $L$ (implying that $Y = X = L$). For let $p \in L$ be a 
boundary point of $Y$, then $\phi$ is a homomorphism when restricted to a 
suitable neighborhood $W$ of $p$, and $W\cap Y \ne \emptyset$ is open in $L$. 
This set is mapped into $M$,
and the same then holds for  $L\cap W$. 
Thus $L \cap W \subseteq Y$, contrary to the choice of $p$.
\epf

\ble\label{two}
Let $\CS$ be a stable plane and suppose that
$U_i$, $i\in\{1,2\}$ are two classical open subsets of $S$, witnessed by homomorphisms 
$\phi_i: U_i \to P_2 \BF$. \\
{\rm a)} There is a local homomorphism $\phi: U_1 \cup U_2 \to P_2 \BF$ that agrees with 
$\phi_1$ on $U_1$. \\
{\rm b)} If a line $L\in \CL$ meets the intersection $V = U_1\cap U_2$, then $L^\phi$ 
is contained in a line.
\ele

Of course, assertion (b) may fail for lines not meeting $U_1 \cap U_2$. 
In other words $\phi$ is, in general, not a homomorphism. 
For a concrete example where this happens see Example 
\ref{extwo}. This example occurs in a Moulton projective 
plane, in fact, in one of its affine subplanes, which is a Salzmann plane as considered in 
\cite{polley-szeb} . In fact, we shall never use assertion (b); instead, we rely on 
Lemma \ref{zush}.  \

Assertion (a) is true and can be proved in the same way if $\phi_1$ and $\phi_2$ are merely 
local homomorphisms, provided that they induce the same map on $U_1\cap U_2$ and moreover
this map is a homomorphism, not just a local one. The last condition is indispensible, as 
shown by Example \ref{extwoplus}. \ 

It would be convenient if one could prove Theorem \ref{main} using Zorn's Lemma. 
The idea would be to take a local homomorphism $U \to P_2\BF$ 
with maximal domain $U$ and to use Lemma \ref{two} in order to prove that $U = S$.
But because of the 
example just mentioned, any attempt to do so is bound to fail.

\bpf
If $V=\emptyset$, then it suffices to define $\phi = \phi_i$ on $U_i$. If $V$ is nonempty, 
then $\phi_2^{-1} \phi_1$ is a homomorphism sending 
$V^\phi_2$ to $V^\phi_1$. According to \cite{locfund}, this map extends to an isomorphism 
$\psi$ of $\CP_2 \BF$. 
Now define $\phi=\phi_1$ on $U_1$ and $\phi = \phi_2\psi$ on $U_2$. These 
definitions agree on $V$ and we have a local homomorphism,
which proves (a). For (b), observe that $\phi$ induces homomorphisms on both 
$U_1$ and $U_2$, hence $(L\cap V)^\phi$ is contained
in a line $M$, and both $(L\cap U_1)^\phi$ and $(L\cap U_2)^\phi$ are contained in the same line.
\epf

\bthm\label{ext}
Let $\CU = (U,\CL_U)$ be a dense open subplane af a stable plane $\CS = (S, \CL)$, 
and $\CP = (P,\CM)$
a compact projective plane of the same dimension.
Then every homomorphism
$\phi : \CU \to \CP$ extends to a (unique) homomorphism $\overline \phi: \CS \to \CP$.
\ethm

\bpf
By definition of a subplane, $\CL_U$ consists of the restrictions $L\cap U$, $L\in \CL$.
Let $a,b$ be two points of $\CU$ and $K = a\vee b \in \CL$,  and consider the 
complement $V = S \setminus K$ 
and the restriction of $\phi$ to the set $V \cap U = U \setminus K$, which is 
still dense in $S$. 
It suffices to extend this restriction to a homomorphism 
on all of $S$. By continuity, this extension will agree with $\phi$ on $K\cap U$. 
For $x \in V$, we set
$x^{\overline \phi} = (a \vee x)^\Phi \wedge (b\vee x)^\Phi$. By construction, 
the map $\overline \phi$ is continuous and injective and extends the restriction of $\phi$. \

Let $L$ be a line meeting $V$. We claim that $L^{\overline \phi}$ 
is contained in $L^\Phi = L^{\overline \Phi}$.
So let $x$ be any point of $L$. By density, there is a sequence $x_n \in V$ converging to $x$.
Furthermore, there exist lines $L_n$ converging to $L$ that contain $x_n$. Since $L$ meets $V$, 
there are points $v_n, w_n \in L_n \cap V$ converging to distinct points 
$v,w \in (L\cap V) \setminus\{x\}$, respectively. It follows that $x_n^\phi \in L_n^\Phi \to 
L^{\overline \Phi}$, whence $x^{\overline \phi} \in L^{\overline \Phi}$. 

Finally, every line $L$ not meeting $V$ is the limit of some sequence of lines $L_n$ 
that do meet $V$. We have shown that $L_n^{\overline \phi}$ is contained in a line. 
Choose sequences $x_n, y_n \in L_n$ converging to distinct points $x,y \in L$. Then 
$L_n^{\overline \phi} = x_n^{\overline \phi} \vee y_n^{\overline \phi}$ converges to 
$K : = x^{\overline \phi} \vee y^{\overline \phi}$. All other points  $z \in L$ are likewise 
limits of sequences in $L_n$,  and the same reasoning shows that their images under 
$\overline \phi$ converge to points of $K$, so by continuity, $L^{\overline \phi} \subseteq K =
L^{\overline \Phi}$, as desired.
\epf

Equality of dimensions is assumed in Theorem \ref{ext} only because our definition 
of homomorphisms requires it.
The proof works equally for homomorphisms in a broader sense, where
the dimension of $\CP$ may exceed that of $\CS$. 

\section{Isotopies of arcs}\label{arcs}

We shall need to know that two arcs with common end points in the point set $S$ of a 
stable plane are isotopic when it is only given that they are homotopic. Without 
further mention, all homotopies and isotopies will be relative endpoints, i.e., 
endpoints are constant throughout the homotopy. Isotopic means that there is a homotopy
$\alpha: I \times I \to S$ such that each map $\alpha_t: s \mapsto \alpha(s,t)$ 
parametrizes an arc in $S$, with $\alpha_0$ and $\alpha_1$ being the given arcs. 
(As customary, $I = [0,1] \subseteq \BR$ denotes the unit interval.) 

If $S$ is a manifold of dimension $d > 2$, then the desired assertion is a general result due to
J. Martin and D. Rolfsen \cite{martin}. Remember that the point set 
of a locally classical plane is indeed a manifold.
This result is false for boundaryless surfaces in general. 
A counterexample on $S = I \times \BS_1$ 
(the annulus, also known as the cylinder) is given by C. D. Feustel \cite{feustel}. 
Therefore we need some
precautions and we take care to give a complete proof, expanding on some brief hints given 
in \cite{polley-szeb} and \cite{polley-moeb}. 
In order to make those arguments work,
a special role will be assigned to \it piecewise $\CS$-linear arcs \rm in a 
2-dimensional stable plane $\CS$. These will be arcs composed of finitely many closed 
intervals contained in lines of $\CS$. We can now state 

\bthm\label{iso}
{\rm \cite{martin}, \cite{polley-szeb}}
Let $\alpha_1$ and $\alpha_2$ be two parametrized 
arcs with common end points in 
the point set $S$ of a stable plane $\CS$. Assume that the point set is a manifold and that 
$\alpha_1, \alpha_2$ are homotopic (rel endpoints). Then each of the following 
conditions implies that the given arcs are isotopic. \

{\rm a)} $S$ is a manifold of dimension $d > 2$, or\

{\rm b)} $S = \BR^2$ and the arcs are piecewise $\CS$-linear.
\ethm

\bpf
The hard part is the proof under condition (a), which is given in \cite{martin}. 
Now consider condition (b). The two arcs may have some intervals of $\CS$-lines in 
common. Apart from that, they have only finitely many intersection points, because 
two $\CS$-lines intersect at most once. Without loss of generality, we may assume 
that the two arcs do not have an initial interval in common. Then the first  point of
$\alpha_1\bigl ( ]0,1] \bigr )$ lying on $\alpha_2 (I)$ is 
$\alpha_1(t_1) = \alpha_2(t_2)$ with $t_1 >0$.
Then the subarcs $\alpha_1[0,t_1]$ and $\alpha_2[0,t_2]$ 
have exactly their end points in common.
In $S = \BR^2$ these arcs together bound a disc by the classical Schoenflies Theorem, 
and it follows that the subarcs are isotopic
rel endpoints. Proceeding inductively, we obtain that $\alpha_1$ and $\alpha_2$ are isotopic. 
\epf

\section{Main results}\label{results}

\bthm\label{main}
Let $\CS = (S,\CL)$ be a locally classical stable plane and assume that
the set of connected lines is dense in $\CL$. 
Then  $\CS$ is classical in 
each of the following cases.
\begin{itemize}
{\rm a)} $\dim S > 2$ and $S$ is simply connected, or\

{\rm b)} $\dim S = 2$ and all lines of $\CS$ are connected.\

{\rm c)} $S$ is either
the Euclidean plane $\BR^2$ or 
the real projective plane ${\bf P} = P_2\BR$.
\end{itemize}
\ethm

Some implications of the hypotheses are obvious: The point set $S$ of a locally 
classical plane is a topological manifold.
For stable planes in general, this is not known.
Density of the set of connected lines implies that $S$ is connected. 
If $S = \bf P$, then the plane is projective and all lines are connected.

The point set of a stable plane satisfying condition (b) is either
$\BR^2$ or  ${\bf P}$, or the Moebius 
band  ${\bf M} = {\bf P} \setminus \{p\}$,  $p \in \bf P$, see
\cite{revisited}, \cite{loe-soyturk}.
The latter possibility is not 
listed in condition (c) because planes on the Moebius band containing disconnected lines
do not always admit compact lines, which are essential to our proof 
(as well as to the proof indicated by Polley \cite{polley-moeb}). So here a natural question 
remains unanswered. 

The point of condition (c) is 
to exclude planes with point set the cylinder $\bf C =\BR^2 \setminus \{0\}$, 
where the theorem is false. See Example  \ref{excylinder}  for a counterexample.  
Section \ref{ex} also contains further examples demonstrating the 
indispensability of most of our hypotheses. \\

\bcor \label{projective}
A locally classical projective plane is classical.
\ecor

\bpf
This follows because projective planes of dimension $>2$ are connected 
and simply connected \cite{CPP} 51.28,   and projective planes have connected lines. 
\epf

\noindent \it Proof of Theorem \ref{main}. \rm
First observe that it suffices to construct a local homomorphism $\phi: \CS \to \CP = \CP_2\BF$.
Indeed, by Lemma \ref{zush}, $\phi$ then maps every connected line of $\CS$ 
into a line of $\CP$.
By density of the set of connected lines, this carries over to all lines, and $\phi$ 
is a homomorphism. (Observe here that every triple of points in a disconnected 
line may be approximated by triples contained in connected lines.)
By Lemma \ref{inj}, $\phi$ is injective and is an isomorphism of $\CS$ onto 
the open subplane $\CS^\phi \le \CP$.\

First we shall prove the theorem assuming that $\dim S > 2$ or $S \approx \BR^2$. Once this 
has been established, the remaining two possibilities $S \approx \bf M$ and $S \approx \bf P$ 
will be treated by reduction to the known case $S \approx \BR^2$. \

Let $A\subseteq S$ be an arc. By compactness of $A$, and because $\CS$ is locally 
classical, there exists a covering of $A$
by finitely many classical connected open subsets $U_1, ... U_k$ of $S$,  
such that $U_i\cap U_j$ is nonempty 
if and only if $\vert i-j\vert\le 1$,
and such that one end point of $A$ is contained in $U_1$ and the other one in $U_k$. 
A covering of $A$ with these properties 
will be called a \it perfect cover\rm . \

At this point, we start constructing a local homomorphism $\CS \to \CP$. 
Choose a point $x_0$ and a 
homomorphism $\phi_0 : U_0 \to \CP$ defined on a neighborhood $U_0$ of $x_0$.
Given any arc $A$ with end point $x_0$, we can use a perfect cover $\CU =\{U_1, ... ,U_k\}$ 
of $A$ such that $U_1 \subseteq U_0$ 
in order to define inductively, using Lemma \ref{two},
a local homomorphism $\phi_\CU$ on the union  $V = \bigcup \CU$ that agrees with $\phi_0$ 
on $V \cap U_0$. We call $\phi_\CU$ an extension of $\phi_0$
along the arc $A$.  
From Lemma \ref{extuniq} we infer that any two extensions of $\phi_0$ along the same arc $A$ 
agree on some neighborhood of $A$, namely, on the connected component of $A$
in the intersection of their domains. In particular, if we extend $\phi_0$ along another arc 
$A'$ with the same endpoints $x_0, x$ and sufficiently close to $A$, then
the values at $x$ taken by the two extensions coincide.\ 

We can now use Theorem \ref{iso} in order to show that any two extensions of $\phi_0$ along arcs
$A_0,A_1$ joining $x_0$ to $x$ take the same value at $x$; in the case $S \approx \BR^2$, we 
take the precaution of comparing only piecewise $\CS$-linear arcs. By assumption, 
the two arcs are homotopic (rel endpoints, as always). The theorem asserts that they are 
isotopic. The intermediate arcs $A_t$ may not be piecewise $\CS$-linear in the 2-dimensional 
case, but nevertheless, by the result of the previous paragraph, 
the value obtained by extension 
along $A_t$ is locally constant as a function of $t$ and hence constant.\

Finally, we define the extension $\phi: \CS \to \CP$ of $\phi_0$ using this construction: 
$x^\phi$ is defined as the value taken at $x$ by the extension of $\phi_0$ along any 
arc with endpoints $x,x_0$. If $S \approx \BR^2$, for this purpose we use only arcs 
that are piecewise $\CS$-linear. Such arcs with preassigned end points always exist;
indeed, because $\CS$ is locally classical, every point $x$ has a neighborhood  $W$
such that $x$ can be joined to every point in $W$ by an arc contained in $L\cap W$, 
where $L$ is a line. The local convexity property just stated actually holds 
in every stable plane.\

It remains to be shown that $\phi$ as constructed above is a local homomorphism. 
Consider any point $x$ in $S$ and a perfect cover of a piecewise $\CS$-linear 
arc $A$ with end points $x_0$ 
and $x$. Let $U_k$ be the last element of 
the perfect cover. Again by local convexity, we may assume that each point of
$U_k$ can be joined to $x$ within $U_k$ by an arc contained in an $\CS$-line. 
Then we can construct the extension of $\phi_0$ at any point of $U_k$
using a prolongation of $A$ by such an arc in $U_k$. 
In the 2-dimensional case, this is again a piecewise $\CS$-linear arc.
The perfect cover $\CU$ and the local homomorphism $\phi_\CU$ defined on its union
are used without change for all such points.
This shows that $\phi$ restricted to $\bigcup \CU$ equals $\phi_\CU$, which by its
construction induces a homomorphism on $U_k$. 
This ends the proof in the cases under consideration.\

If $S \approx \bf P$, then $S$ is compact and $\CS$ is a projective plane. 
Delete one line of $\CS$ to obtain an affine plane with point set $\BR^2$. 
This plane is classical
by what we already know, and its projective extension $\CS$ is classical, as well.\

The last missing case is that $S$ is a Moebius band,
$S = {\bf M} = {\bf P}\setminus \{p\}$. In this case, we shall 
need to use compact lines of
$\CS$. It is easy to produce examples of stable planes with $S = \bf M$ where neither compact
lines nor connected lines occur. It is not known whether compact lines can be absent 
when the connected lines form a dense set. Therefore, this case will be treated under 
assumption (b), that is, all lines are assumed to be connected.\

If this is assumed, then it is known that every point is incident 
with both compact and non-compact lines (homeomorphic to $\CS_1$ and to $\BR$, 
respectively), and that the set of compact lines is open in the line space; moreover, a 
compact line intersects every other line, see 
\cite{sz-pjm}, \cite{betten}. \

Choose a non-compact line $L$. Like every line,  $L$ is closed in $\bf M$. 
Hence in the one-point compactification ${\bf P} = {\bf M} \cup \{p\}$ 
of $S =\bf M$, the closure of $L$ is a circle $L \cup\{p\}$. The inverse image $\tilde L$ 
of this circle under the two-sheeted covering ${\BS_2}\to \bf P$ is either again a circle or a 
union of two disjoint circles. In the first case, $\tilde L$ contains the two inverse 
images of $p$, and in the second case each circle contains one of the those two points. 
In the second case, one sees that the complement of $\tilde L$ consists of three connected 
components, one cylinder and two disks. It follows that $S \setminus L$ is disconnected (a 
union of a disk and a Moebius band). This will shortly be used to rule out this possibility. 
In the first case, it follows from the Schoenflies Theorem that ${\BS_2} \setminus \tilde L$ is 
the union of two disjoint open disks, each of which must be mapped homeomorphically onto 
$S \setminus L$ by the covering map. Thus, $S \setminus L$ is homeomorphic to $\BR^2$. \

In order to see that $ S \setminus L$ is connected and hence homeomorphic to $\BR^2$, we use 
the compact lines. Take two points $x,y$ not on $L$ and choose compact lines $K_x$ and $K_y$
containing $x$ and $y$, respectively. Then $K_x$ intersects $K_y$, and since compact 
lines form an open set, we may assume that the intersection point does not belong to $L$. 
Thus we have a connected subset $(K_x \cup K_y) \setminus L$ of $S \setminus L$ containing \
both $x$ and $y$. 

On the subplane $\CS_L$ induced on $S_L = S \setminus L \approx \BR^2$, we may use the 
known procedure to construct a local homomorphism $\phi_L: \CS_L \to \CP$
by extending a homomorphism $\phi_0$ defined on some open subset $U_0$ of $S_L$; this does 
not require any connectedness assumptions for lines. However, $\CS_L$ does have too many 
disconnected lines, so we cannot conclude directly that $\phi$ is a 
homomorphism. At least, compact lines of $\CS$ induce 
non-compact connected lines of $\CS_L$, and these are mapped into lines of $\CP$ 
according to Lemma \ref {zush}.\

Now choose a second noncompact line $G$ not meeting $L$. 
Define $\phi_G$ on $S_G$ as before, using the same 
set $U_0 \subseteq S_L \cap S_G$ and the same initial homomorphism $\phi_0: U_0 \to P$. 
We claim that  $\phi_L$ agrees with $\phi_G$ on the intersection 
$S_L \cap S_G = S \setminus \{L \cup G\}$. Thus we obtain a local homomorphism on all of 
$S = S_L \cup S_G$. 
Then we can once more use the hypothesis that all lines of $\CS$ are 
connected and apply Lemma \ref{zush}
in order to see that we have a homomorphism.

The claim above is proved using the compact lines. Every point of $U_0$ belongs to 
some compact line $K$. Since compact lines form an open set, 
it follows that each point of some neighborhood $V$ of $K$ is the intersection of 
two compact lines meeting $U_0$. As $\phi_L$ and $\phi_G$ agree on $U_0$ and 
behave like homomorphisms on compact lines (Lemma \ref{zush}), $\phi_L$ and $\phi_G$ agree 
on $V$. We repeat this argument: every point $x$ is the intersection of two compact lines.
These intersect $K$ and, hence, $V$. Thus, $x^{\phi_L} = x^{\phi_G}$. This completes the proof. 
\ok\bigskip\noindent

\section{Examples}\label{ex}

We shall give a sequence of examples in order to show that our results cannot be 
generalized very much. Some open questions nevertheless remain. 
Almost all examples in this section are open subplanes of the Moulton 
projective planes $\CM_k$, $k>1$.
These are non-classical 2-dimensional planes whose full group of automorphisms $\Sigma_k$ is 
4-dimensional, locally isomorphic to $\mathop{\rm GL}_2\BR$, and of Lenz type III. No other 
compact connected projective planes of Lenz 
type III exist. This was proved by H. Salzmann in each dimension separately;  
see \cite{CPP}, 64.18 for a unified proof. In particular, there is no higher-dimensional 
analogue of the Moulton planes.
There are two useful ways to describe $\CM_k$: \

The original \it Moulton
model \rm \cite{moulton}
is obtained from the real affine plane by replacing lines of negative slope $s < 0$ 
by kinked lines that have slope $s$ in the right half plane $x \ge 0$ and slope $ks$ in the 
left half plane $x \le 0$; the projective closure of the resulting affine plane is the Moulton 
plane $\CM_k$. The automorphism group fixes the $y$-axis $Y$ and the point at infinity $q$ 
on the $x$-axis and is transitive on the remainder $C_k$ of the point set, which is a 
topological cylinder. 

The second, \it radial  model \rm of $\CM_k$ is designed to make the action of $\Sigma_k$ on 
$C_k$ visible. In this model, $Y$ becomes the line at infinity and $q$ is the origin, so 
$C_k$ is just $\BR^2 \setminus \{0\}$ in this model. 
Lines through the origin are ordinary lines of the real affine plane. For a 
complete description, see \cite{CPP}, 31.25 (b) for the Moulton model and Section 34 for the 
radial one; compare also the article \cite{laguerre}, 
which by the way corrects some minor errors
contained in Section 34 of \cite{CPP}. \\

\begin{Example}\label{exdisconn} 
\rm The affine part of the Moulton model is divided by the $y$-axis $Y$,
which is the locus of the kinks, into two open half planes, which are obviously classical. 
Their union is, however, not classical. This can be seen either by applying Theorem \ref{ext}
(a homomorphism to the real projective plane $\CP$ would extend to an 
isomorphism $\CM_k \to \CP$, which is impossible),
or by examining Desargues configurations. So the plane induced on the complement of $Y$ in the 
affine plane is locally classical, but not classical. This shows that connectedness of the 
point set is essential for the results discussed here. 
\end{Example}

\begin{Example}\label{exhilbert}
\rm Hilbert's first example of a non-desarguesian projective plane 
\cite{hilbert}, Kap. V, \S 23, compare also \cite{stroppel}, (14), is obtained 
from the real projective plane $\CP$ by changing the lines inside the unit circle $K$ 
of $\BR^2$. In fact, the modified interior of $K$ is isomorphic to an open subplane of 
$\CP$ under a Moebius transformation. Thus 
as before, the complement of $K$ in Hilbert's plane is locally classical but not classical.
In higher dimensions, no analogue of Hilbert's plane is known.
\end{Example}

\begin{Example}\label{excylinder} 
\rm The open subplane $\CC_k$ of $\CM_k$ induced on the cylinder $C_k$ 
(the complement of the two fixed elements of $\Sigma_k$) is connected, and the connected lines 
form a dense subset (only the lines passing through the fixed point $q$ 
become disconnected). The plane $\CC_k$ is 
locally classical, because the (classical) right half plane of the Moulton model 
corresponds in $\CC_k$ to one half of the cylinder, bounded by some ordinary line through $q$. 
The action of the automorphism group shows that every half of $\CC_k$ defined by such a line 
is classical.  So $\CC_k$ is locally classical, but not classical for the usual reasons.
This example shows that Theorem \ref{main} is false for planes with point 
set the cylinder (necessarily containing disconnected lines). 
Yet it is true if the point set is the Moebius band and lines are connected. Both these 
surfaces are connected but not simply connected, and the narrow distinction between the 
two seemed another good reason to thouroughly examine the existing proofs of Theorem \ref{main}. 
\end{Example}

\begin{Example}\label{extwo}
\rm Like in the previous example, start with one half $U$ of $C_k$ 
bounded by a line through $q$. Rotate this half plane though an angle of 90 degrees. This can 
be done by an automorphism of $\CM_k$, so we obtain another classical open subset $V$. 
The union $U \cup V$, however, does not define a classical subplane. This is best seen in 
the Moulton model, where $U$ corresponds to the right open half plane and $V$ 
to the set of points $(x,y)$ with $xy >0$ plus the positive part of the line at infinity. 
The union $U \cup V$ contains non-closing Desargues configurations
and hence is not classical. \

Applying Lemma \ref{two}, we obtain a local homomorphism $\phi$ of $U \cup V$ 
into the real projective plane
which cannot be a homomorphism. This shows in particular that in applications of part (b) of 
this Lemma, the hypothesis that the line 
in question meets the intersection of the given open sets cannot be ignored. 
Note that $U\cup V$ is contained in the affine subplane of $\CM_k$ obtained by deleting the 
line $Y$ fixed by all automorphisms. This is a Salzmann 
plane as considered in \cite{polley-szeb}. 
\end{Example}

\begin{Example}\label{extwoplus}
\rm Let $W$ be another rotated half of $C_k$ and $\psi$ a homomorphism of $W$ into $\CP_2\BR$.
From Example \ref{extwo}, import the set $U \cup V$ and the local 
homomorphism $\phi$ defined on it. 
We may choose $W$ such that $U \cup V \cup W$ is all of $C_k$. We try to apply the sharpened
version of Lemma \ref{two} to the maps $\phi$ and $\psi$. 
If this Lemma would apply, we would get a local 
homomorphim $\CC_k \to \CP_2\BR$. This would in fact be a homomorphism, 
because the set of connected lines is dense. 
That is impossible, so we see that the sharpened Lemma does 
need the hypothesis that the two given local homomorphism induce a true homomorphism on the 
intersection of their domains.
\end{Example}

\begin{Example}\label{exsimplyconn}
\rm As in Example \ref{extwo}, start with $\CC_k$ and delete only 
a half line joining 
$q$ to the line at infinity. The remainder $R$ is homeomorphic to $\BR^2$, in particular, 
simply connected, and it defines a 
subplane which is locally classical but not classical. In this example, there is an 
open set of disconnected lines. 
\end{Example}

\begin{Example} \label{exspiral}
\rm Another example with the same properties as Example \ref{exsimplyconn} can be 
obtained by deleting from $C_k$ an orbit of one of the one-parameter group of 
\it spiral rotations \rm sending $re^{i\phi}$ to $e^{st} e^{i(\phi + t)}$, 
see \cite{CPP}, 34.4. 
The parameter $s$ here is related to $k$ by $k = e^{2\pi s}$.  This example has the 
extra benefit that it admits the one-dimensional group of 
spiral rotations. 
\end{Example}

The remaining \it open questions \rm concern the following situations: (1) $\dim S> 2$ and $S$ 
is not simply connected, and (2) $S = \bf M$ and there are disconnected lines. 
Our examples do not exclude the possibility that Theorem \ref{main} remains true 
in these cases, but we have no positive result. Moreover, in 
higher dimensions there is a total lack of counterexamples of the type discussed here. 
\\

\bibliographystyle{plain}

\bigskip
\bigskip
\noindent{Rainer L\"owen\\ 
Institut f\"ur Analysis und Algebra\\
Technische Universit\"at Braunschweig\\
Universit\"atsplatz 2\\
38106 Braunschweig\\
Germany}

\end{document}